\documentclass[12pt]{article}

\usepackage{fullpage}
\usepackage{amssymb}
\usepackage{amsmath}
\usepackage{amsthm}
\usepackage[usenames]{color}
\usepackage{graphicx}
\usepackage{booktabs}
\usepackage{float}
\usepackage{enumitem} 
\usepackage{longtable} 
\usepackage[colorlinks=true,
linkcolor=webgreen,
filecolor=webbrown,
citecolor=webgreen]{hyperref}

\definecolor{webgreen}{rgb}{0,.5,0}
\definecolor{webbrown}{rgb}{.6,0,0}

\theoremstyle{plain}
\newtheorem{theorem}{Theorem}

\newtheorem{lemma}[theorem]{Lemma}
\newtheorem{proposition}[theorem]{Proposition}

\theoremstyle{definition}
\newtheorem{definition}[theorem]{Definition}

\newtheorem{conjecture}[theorem]{Conjecture}

\theoremstyle{remark}
\newtheorem{remark}[theorem]{Remark}

\title{An effective analytic recurrence for prime numbers}

\author{Benoit Cloitre\\
Email: \texttt{benoit.cloitre@proton.me}}

\date{}

\begin{document}

\maketitle

\begin{abstract}
The Golomb--Keller formula expresses the next prime $p_{n+1}$ as a recurrence relation in terms of the first $n$ primes $p_1, \ldots, p_n$ using the Riemann zeta function and an Euler product, but requires taking a limit as $s \to \infty$, rendering it non-constructive. We transform this asymptotic formula into an effective recurrence by proving that a finite parameter $s \leq p_n$ suffices when combined with the ceiling function, establishing a constructive method valid for all $n \geq 1$.

The minimal integer parameter $s_n$ (OEIS A389650) reveals deep connections to prime constellations. We prove $\liminf_{n\to\infty} \sigma_n = 0$ unconditionally, where $\sigma_n = s_n/p_n$. The limit superior $C = \limsup \sigma_n$ satisfies $\log \psi \lesssim C \leq 0.4332$, where $\psi \approx 1.46557$ is the supergolden ratio. The lower bound is conditional on the twin prime conjecture; the upper bound is unconditional. The constant $C$ encodes the densest admissible prime constellation occurring infinitely often, linking our formula to the Hardy--Littlewood conjectures.

The method extends to Dirichlet L-functions, yielding other effective formulas for calculating $p_{n+1}$ but also for predicting residues of $p_{n+1}$ modulo any integer with reduced precision requirements.
\end{abstract}



\section{Introduction}

Finding explicit formulas for primes is a classical problem in number theory. While direct, non-recursive formulas exist (e.g., Willans, Mills), they are computationally infeasible \cite{Ribenboim}. We focus instead on recurrence relations expressing $p_{n+1}$ in terms of $p_1, \ldots, p_n$ that, while also computationally intensive, reveal structural properties of the prime sequence.

Gandhi \cite{Gandhi1971} provided the first such recurrence using primorials and the Möbius function. Vanden Eynden \cite{VandenEynden1972} simplified the proof, Jakimczuk \cite{Jakimczuk2024} generalized the approach, and Trefeu \cite{Trefeu2025} clarified the link to Eratosthenes' sieve using generating functions.

Golomb \cite{Golomb1976} used analytic number theory, independently rediscovered by Keller \cite{Keller2007}, to derive:
\[
p_{n+1} = \lim_{s\to\infty} \left[\left(\prod_{k=1}^n \left(1-\frac{1}{p_k^s}\right)\right)\zeta(s)-1\right]^{-1/s}.
\]
This formula is impractical due to the limit.

\subsection{The effective approach}

Keller's derivation truncated the sum $\sum_{j=1}^{2p_n-1} j^{-s}$ in place of $\zeta(s)$ while taking $s \to \infty$. We take the opposite approach: retain the full series $\zeta(s)$ computed to working precision, but use finite $s$. This truncates the exponent rather than the sum, making the formula constructive. The complete series $\zeta(s)$ provides the global structure needed to bound the error, while we prove that finite $s \leq p_n$ suffices to isolate $p_{n+1}$.

\subsection{Main results}

For each $n \geq 1$, let $s_n$ denote the smallest integer satisfying
\[
p_{n+1} = \left\lceil \left( -1 + \zeta(s_n) \prod_{j=1}^{n} \left(1 - \frac{1}{p_j^{s_n}}\right) \right)^{-1/s_n} \right\rceil.
\]
Such an integer exists (Proposition \ref{prop:existence_sn}). We establish:

\begin{enumerate}
\item $s_n \leq 2p_n$ for all $n \geq 1$, using Bertrand's postulate (Theorem \ref{thm:bertrand}).
\item $s_n \leq p_n$ for all $n \geq 1$, using Nagura's theorem (Theorem \ref{thm:nagura}).
\end{enumerate}

The sequence $(s_n)_{n \geq 1}$ is recorded as sequence A389650 in the OEIS \cite{OEIS}. Writing $\sigma_n = s_n/p_n$, we prove:

\begin{enumerate}
\setcounter{enumi}{2}
\item $\liminf_{n \to \infty} \sigma_n = 0$ unconditionally (Proposition \ref{prop:liminf}).
\item $C := \limsup_{n \to \infty} \sigma_n$ satisfies $0.3823 \lesssim C \leq 0.4332$. The lower bound, conditional on the twin prime conjecture, approximately equals $\log \psi \approx 0.3823$. The upper bound is unconditional (Theorems \ref{thm:twin}, \ref{thm:limsup}).
\item For any fixed $c > c_0 \approx 0.5956$, the formula with $s = cp_n$ holds for all sufficiently large $n$ (Theorem \ref{thm:fixed}). We conjecture this holds for any $c > C$ (Conjecture \ref{conj:general}).
\item The method extends to Dirichlet L-functions, allowing prediction of residues $p_{n+1} \pmod 4$ with reduced precision (Theorems \ref{thm:dirichlet_main}, \ref{thm:dirichlet_sign}).
\end{enumerate}

Beyond computational applications, the sequence $(s_n)$ reveals intrinsic number-theoretic structure through its connection to prime constellations, gap distributions, and the Hardy--Littlewood conjectures. Our data (Appendix A, Figures \ref{fig:ratio}--\ref{fig:histogram_trimmed}) illustrate the behavior of $\sigma_n$, its correlation with specific gap configurations, and its empirical distribution.

\section{Golomb's unifying framework}

Golomb \cite{Golomb1976} showed that various prime recurrence formulas fit a single framework involving a probability distribution $\alpha(k)$ on $\mathbb{N}$ and an operator $T$.

Let $P_n = \prod_{j=1}^{n} p_j$ be the $n$th primorial. Define
\[
\gamma(P_n) = \sum_{\substack{k=1 \\ \gcd(k,P_n)=1}}^\infty \alpha(k).
\]

Taking $\alpha(k) = (b-1)b^{-k}$ (geometric distribution) and $T(x) = \lfloor -\log_b(x) \rfloor + 1$ recovers Gandhi's formula (when $b=2$).

Taking $\alpha(k) = k^{-s}/\zeta(s)$ for $s>1$, the Euler product gives
\[
\gamma(P_n) = \frac{1}{\zeta(s)} \sum_{\gcd(k,P_n)=1} k^{-s} = \prod_{j=1}^n \left(1 - \frac{1}{p_j^s}\right).
\]
The operator $T(x) = \lim_{s\to\infty} x^{-1/s}$ extracts the leading index $p_{n+1}$, yielding the Golomb--Keller formula.

\section{From asymptotic to effective formula}

\subsection{Analytical tools}

\begin{lemma}[Integral bounds]
\label{lem:integral}
Let $s > 1$ and $m \geq 2$ be an integer.
\begin{enumerate}
\item $\displaystyle\sum_{k=m}^{\infty} k^{-s} < \frac{(m-1)^{1-s}}{s-1}$.
\item $\displaystyle\sum_{k=m}^{\infty} k^{-s} \le m^{-s} \left(1 + \frac{m}{s-1}\right)$.
\end{enumerate}
\end{lemma}

\begin{proof}
(1) follows from standard integral comparison: $\sum_{k=m}^\infty k^{-s} < \int_{m-1}^\infty x^{-s}\,dx$.
(2) follows by isolating the first term: $m^{-s} + \sum_{k=m+1}^\infty k^{-s} \leq m^{-s} + \int_m^\infty x^{-s}\,dx = m^{-s} + m^{1-s}/(s-1)$.
\end{proof}

\subsection{Dirichlet series machinery}

\begin{definition}
For $\Re(s) > 1$, define
\begin{align*}
D_n(s) &= \sum_{\substack{k \geq 1 \\ \gcd(k,P_n)=1}} k^{-s}, \\
T_n(s) &= \sum_{\substack{k \geq p_{n+2} \\ \gcd(k,P_n)=1}} k^{-s}, \\
T_n^*(s) &= p_{n+1}^s T_n(s).
\end{align*}
Observe that $D_n(s) - 1 = p_{n+1}^{-s} + T_n(s)$. The condition $\gcd(k, P_n)=1$ implies that $k$ is either 1 or composed only of primes exceeding $p_n$.
\end{definition}

\begin{lemma}[Euler product]
\label{lem:euler}
For $\Re(s) > 1$,
\[
D_n(s) = \zeta(s) \prod_{j=1}^{n} (1 - p_j^{-s}).
\]
\end{lemma}

\begin{proof}
The Euler product for $\zeta(s)$ (see \cite{Apostol}, Theorem 11.6) is $\prod_{p} (1-p^{-s})^{-1}$. Splitting this product:
\[
\zeta(s) = \left(\prod_{j=1}^{n}(1-p_j^{-s})^{-1}\right) \cdot \left(\prod_{j=n+1}^{\infty}(1-p_j^{-s})^{-1}\right).
\]
The second infinite product equals $D_n(s)$ as it represents the sum over integers whose prime factors exceed $p_n$. Multiplying by $\prod_{j=1}^{n} (1-p_j^{-s})$ yields the result.
\end{proof}

\subsection{Convergence properties}

Define $h(s) = (D_n(s)-1)^{-1/s}$.

\begin{lemma}
\label{lem:convergence}
For any $n \ge 1$:
\begin{enumerate}
\item $h(s) < p_{n+1}$ for all $s > 1$.
\item $\lim_{s\to\infty} h(s) = p_{n+1}$.
\end{enumerate}
\end{lemma}

\begin{proof}
(1) Since $T_n(s)>0$, $D_n(s) - 1 > p_{n+1}^{-s}$. Raising to the power $-1/s$ reverses the inequality.

(2) Write $h(s) = p_{n+1} (1 + T_n^*(s))^{-1/s}$. We show $T_n^*(s) \to 0$ as $s \to \infty$. Using Lemma \ref{lem:integral}(1),
\[
T_n^*(s) < \frac{p_{n+2}-1}{s-1} \left(\frac{p_{n+1}}{p_{n+2}-1}\right)^s.
\]
We verify $p_{n+1} < p_{n+2} - 1$ for $n \geq 1$. For $n=1$, $3 < 5-1$. For $n \geq 2$, $p_{n+1}$ and $p_{n+2}$ are odd, so $p_{n+2} - p_{n+1} \geq 2$, giving $p_{n+2}-1 \geq p_{n+1}+1$.

Since the ratio $p_{n+1}/(p_{n+2}-1)$ is strictly less than 1, the bound decays exponentially in $s$. Thus $T_n^*(s) \to 0$ and $h(s) \to p_{n+1}$.
\end{proof}

\begin{lemma}[Monotonicity of $h$]
\label{lem:h_monotone}
The function $h(s)$ is strictly increasing on $(1, \infty)$.
\end{lemma}

\begin{proof}
Let $F(s)=D_n(s)-1$. Analyzing the derivative, $h'(s)/h(s) = \frac{1}{s^2}(\log F(s) - s F'(s)/F(s))$. Thus $h'(s) > 0$ is equivalent to $\log F(s) > s F'(s)/F(s)$.

Writing $F(s) = p_{n+1}^{-s}(1 + T_n^*(s))$, the logarithmic derivative yields the equivalent condition
\begin{equation}
\label{eq:monotonicity_condition}
\log(1 + T_n^*(s)) > s \frac{(T_n^*)'(s)}{1 + T_n^*(s)}.
\end{equation}
The left side is positive since $T_n^*(s)>0$. The derivative $(T_n^*)'(s)$ is
\[
(T_n^*)'(s) = \sum_{\substack{k \geq p_{n+2} \\ \gcd(k,P_n)=1}} \left(\frac{p_{n+1}}{k}\right)^s \log\left(\frac{p_{n+1}}{k}\right).
\]
Since $k \geq p_{n+2} > p_{n+1}$, the logarithms are negative. Thus $(T_n^*)'(s) < 0$. The right side is negative, so the inequality holds.
\end{proof}

\begin{proposition}[Existence and characterization of $s_n^*$]
\label{prop:existence_sn}
For each $n \geq 1$, there exists a unique $s_n^* > 1$ such that $h(s_n^*) = p_{n+1} - 1$. Furthermore, $h(s) > p_{n+1} - 1$ if and only if $s > s_n^*$.
\end{proposition}

\begin{proof}
By Lemma \ref{lem:convergence}, $\lim_{s\to\infty} h(s) = p_{n+1}$. As $s \to 1^+$, $D_n(s) \to \infty$ (by Lemma \ref{lem:euler}), so $h(s) \to 0$. Since $h$ is continuous and strictly increasing (Lemma \ref{lem:h_monotone}), the Intermediate Value Theorem guarantees a unique $s_n^*$ such that $h(s_n^*) = p_{n+1}-1$.
\end{proof}

\begin{lemma}[Connection to the ceiling formula]
\label{lem:sn_ceiling}
For any $s > 1$, $\lceil h(s) \rceil = p_{n+1}$ if and only if $s > s_n^*$.
The smallest integer satisfying this condition is $s_n := \lfloor s_n^* \rfloor + 1$.
\end{lemma}

\begin{proof}
$\lceil h(s) \rceil = p_{n+1}$ is equivalent to $p_{n+1} - 1 < h(s) \leq p_{n+1}$. Since $h(s) < p_{n+1}$, this reduces to $h(s) > p_{n+1} - 1$. By Proposition \ref{prop:existence_sn}, this holds if and only if $s > s_n^*$.
\end{proof}

\begin{definition}[Minimal parameter values]
\label{def:minimal}
Define $s_n = \lfloor s_n^* \rfloor + 1$, $\sigma_n = s_n/p_n$, and $\sigma_n^* = s_n^*/p_n$.
\end{definition}

\begin{lemma}[Asymptotic equivalence]
\label{lem:equiv}
We have $0 < \sigma_n - \sigma_n^* < 1/p_n$. Thus $\sigma_n - \sigma_n^* \to 0$ as $n \to \infty$, implying their liminf and limsup coincide.
\end{lemma}

\subsection{First effective bound}

\begin{theorem}[Bertrand's postulate bound]
\label{thm:bertrand}
For all $n \ge 1$, we have $s_n \leq 2p_n$.
\end{theorem}

\begin{proof}
We show that for $s = 2p_n$, $h(s) > p_{n+1} - 1$, equivalent to $D_n(s)-1 < (p_{n+1}-1)^{-s}$.
We bound $D_n(s)-1$ using Lemma \ref{lem:integral}(1) with $m=p_{n+1}$:
\[
D_n(s) - 1 < \sum_{k=p_{n+1}}^{\infty} k^{-s} < \frac{(p_{n+1}-1)^{1-s}}{s-1}.
\]
The desired inequality holds if $(p_{n+1}-1)^{1-s}/(s-1) < (p_{n+1}-1)^{-s}$. This simplifies to $p_{n+1}-1 < s-1$, or $s > p_{n+1}$.

Bertrand's postulate (see \cite{HardyWright}, Theorem 418) asserts $p_{n+1} < 2p_n$ for all $n \geq 1$. Choosing $s = 2p_n$ satisfies $s > p_{n+1}$. This proves $h(2p_n) > p_{n+1}-1$. By Proposition~\ref{prop:existence_sn} and the strict monotonicity of $h$ (Lemma~\ref{lem:h_monotone}), this implies $2p_n > s_n^*$. Since $s_n = \lfloor s_n^* \rfloor + 1$ and $2p_n$ is an integer, $s_n^* < 2p_n$ implies $\lfloor s_n^* \rfloor \leq 2p_n - 1$, giving $s_n \leq 2p_n$.
\end{proof}

\subsection{Worked example: computing $p_7$ from $n=6$}

We compute $p_7 = 17$ starting from $n=6$ (primes up to 13). We use the Euler product (Lemma \ref{lem:euler}) to compute $D_n(s)$.

\begin{remark}
\label{rem:precision_example}
When $s$ is large, both $\zeta(s)$ and the product approach 1, leading to catastrophic cancellation when computing $D_n(s)-1$. This requires high working precision. The computations in this paper rely on high-precision arithmetic libraries such as PARI/GP \cite{PARI2023} or Python's \texttt{mpmath}. For $n \leq 200$ (Appendix A), 100 digits suffice. However, for exploring larger $n$ (e.g., $n \approx 500$), approximately 2500 digits of working precision are required due to increasing cancellation as $s$ grows.
\end{remark}

Using $s = 2p_6 = 26$ (with high precision):
\begin{align*}
D_6(26) - 1 &\approx 1.11344499506 \times 10^{-32}, \\
h(26) &\approx 16.941817904\ldots, \quad p_7 = \lceil h(26) \rceil = 17.
\end{align*}

The actual minimal value is $s_6 = 8$ (see Appendix A), giving $h(8) \approx 16.5189076\ldots$.

\section{Sharp bound and asymptotic analysis}

\subsection{Unconditional proof of $s_n \le p_n$}

We use Nagura's theorem \cite{Nagura1952}: for $x \geq 25$, there is a prime in $(x, 1.2x)$.

\begin{theorem}[Sharp bound]
\label{thm:nagura}
For all $n \geq 1$, we have $s_n \leq p_n$.
\end{theorem}

\begin{proof}
For $n \leq 9$, direct numerical verification (see Appendix A) confirms $s_n \leq p_n$.

For $n \geq 10$ ($p_n \geq 29$), Nagura's theorem applies. Let $s = p_n$. We must show $T_n(s) < (p_{n+1}-1)^{-s} - p_{n+1}^{-s}$. Dividing by $p_{n+1}^{-s}$, we verify $R(n) < \text{RHS}(n)$, where
\begin{align*}
R(n) &= \left(\frac{p_{n+1}}{p_{n+2}}\right)^{p_n} \left(1 + \frac{p_{n+2}}{p_n-1}\right), \\
\text{RHS}(n) &= \left(1 - \frac{1}{p_{n+1}}\right)^{-p_n} - 1.
\end{align*}
We used the integral bound (Lemma \ref{lem:integral}(2)) to bound $T_n^*(s)$ by $R(n)$.

Write $R(n) = A_n \cdot B_n$. Since $n \geq 10$, $g_{n+1} = p_{n+2} - p_{n+1} \geq 2$.
By Nagura, $p_{n+2} < 1.2 p_{n+1} < 1.44 p_n$.
\[
\frac{p_{n+1}}{p_{n+2}} \leq 1 - \frac{2}{p_{n+2}} < 1 - \frac{2}{1.44 p_n} = 1 - \frac{25}{18p_n}.
\]
Using $(1-x)^n \leq e^{-nx}$,
\[
A_n \leq \exp\left(-\frac{25}{18}\right) \approx 0.2493 < 0.25.
\]

For $B_n = 1 + p_{n+2}/(p_n-1)$. Since $p_n \geq 29$, we have $p_n/(p_n-1) \leq 29/28$. Thus,
\[
B_n < 1 + 1.44 \frac{p_n}{p_n-1} \leq 1 + \frac{144}{100} \times \frac{29}{28} = 1 + \frac{36}{25} \times \frac{29}{28} = 1 + \frac{261}{175} \approx 2.4914 < 2.5.
\]
Combining gives $R(n) < 0.25 \times 2.5 = 0.625$.

For the right side, using $p_{n+1} < 1.2p_n$ and $-\log(1-x) > x$ for $x \in (0,1)$,
\begin{align*}
\text{RHS}(n) &> \exp\left(\frac{p_n}{p_{n+1}}\right) - 1 > e^{1/1.2} - 1 = e^{5/6} - 1 \approx 1.3016.
\end{align*}
Since $0.625 < 1.30$, the inequality holds for all $n \geq 10$.
\end{proof}

\subsection{Asymptotic behavior of $\sigma_n$}

We now analyze the asymptotic behavior of $\sigma_n$, establishing unconditional bounds and conditional results related to prime constellations.

\begin{proposition}
\label{prop:liminf}
We have $\liminf_{n \to \infty} \sigma_n = 0$.
\end{proposition}

\begin{proof}
By Lemma \ref{lem:equiv}, it suffices to prove $\liminf \sigma_n^* = 0$. Suppose for contradiction that $\liminf \sigma_n^* = \epsilon > 0$. Then for large $n$, $s_n^* \geq (\epsilon/2) p_n$.

By definition, $1 + T_n^*(s_n^*) = (1 - 1/p_{n+1})^{-s_n^*}$.
Using $p_n/p_{n+1} \geq 1/2$ for large $n$ (by PNT or Bertrand's postulate) and $(1-x)^{-c} > e^{cx}$ for $x \in (0,1), c>0$,
\[
1 + T_n^*(s_n^*) > \exp\left(\frac{s_n^*}{p_{n+1}}\right) \geq \exp\left(\frac{\epsilon}{4}\right).
\]
Thus, $T_n^*(s_n^*) \geq e^{\epsilon/4} - 1 > 0$.

We derive an upper bound using Lemma \ref{lem:integral}(1). Let $G = g_{n+1}$.
\[
T_n^*(s) < \frac{p_{n+2}-1}{s-1} \left(\frac{p_{n+1}}{p_{n+1}+G-1}\right)^s.
\]
For large $n$, $s_n^*-1 \geq (\epsilon/4)p_n$. Since $G = o(p_n)$ by PNT, the coefficient $(p_{n+2}-1)/(s_n^*-1)$ is bounded by $K_\epsilon = 5/\epsilon$.

The exponential term is bounded by $\exp(-\frac{\epsilon}{4}(G-1))$.
Combining gives $T_n^*(s_n^*) < K_\epsilon \exp(-\frac{\epsilon}{4}(g_{n+1}-1))$.

The lower and upper bounds imply $g_{n+1}$ is bounded by a constant depending only on $\epsilon$. This contradicts Westzynthius's theorem \cite{Westzynthius1931} ($\limsup g_n = \infty$).
\end{proof}

\begin{theorem}[Limsup bound]
\label{thm:limsup}
Let $C = \limsup_{n \to \infty} \sigma_n$. Then unconditionally, $C \leq 0.4332$.
\end{theorem}

\begin{proof}
Fix $c = 0.4332$ and $s = cp_n$. We use Dusart's explicit bounds \cite{Dusart2010} (Theorem 6.9). For $p_n \geq X_0 = 396738$, $p_{n+1}/p_n \leq 1/a_0$, where $a_0 = 1 - 1/(25\log^2 X_0) \approx 0.999760$.

We decompose $T_n^*(s) = P_n^*(s) + C_n^*(s)$ (prime and composite terms). (Note: $P_n^*(s)$ denotes the normalized prime tail sum and is unrelated to the primorial $P_n$.)

The smallest composite term is $p_{n+1}^2$. Using the integral bound, one verifies that $C_n^*(s)$ is exceedingly small ($< 10^{-956994}$ for $p_n \geq X_0$).

We analyze the prime terms $P_n^*(s)$. For primes $p > 3$, admissibility modulo 3 dictates the densest pattern of gaps. A sequence of consecutive gaps cannot be $(2, 2)$, as one of $p, p+2, p+4$ must be divisible by 3 (unless $p=3$). The densest admissible pattern is therefore $(2, 4, 2, 4, \ldots)$. This follows from the admissibility constraints for prime $k$-tuples (see \cite{HardyLittlewood1923}, \S3).

Let $S_j$ be the cumulative gaps. This pattern gives $S_{2k-1} = 6k - 4$ and $S_{2k} = 6k$.

Let $x_j = S_j/p_{n+1}$. We use $(1 + x_j)^{-s} \leq \exp(-sx_j/(1+x_j))$.
We have $s/p_{n+1} = c p_n/p_{n+1} \geq c a_0$. We analyze the factor $1/(1+x_j)$. The sum $P_n^*(s)$ is dominated by small $j$. We introduce a cutoff $y_0 = 10^{-2}$. For terms where $S_j/p_{n+1} \leq y_0$, we have $1/(1+x_j) \geq 1/(1+y_0)$.

Thus, for these dominant terms, $s x_j/(1+x_j) \geq \kappa c S_j$, where $\kappa = a_0/(1+y_0) \approx 0.989861$. The contribution from the tail where $x_j > y_0$ is negligible due to the rapid exponential decay: $e^{-s x_j} < e^{-c a_0 y_0 p_{n+1}}$. Since $p_n \geq X_0$, this tail is bounded by $e^{-1718} < 10^{-746}$. We therefore approximate $P_n^*(s)$ using $\kappa$ for the entire sum.

\begin{align*}
P_n^*(s) &\lesssim \sum_{k=1}^{\infty} e^{-\kappa c(6k-4)} + \sum_{k=1}^{\infty} e^{-\kappa c \cdot 6k} \\
&= (e^{4\kappa c} + 1) \cdot \frac{e^{-6\kappa c}}{1 - e^{-6\kappa c}} = \frac{e^{-2\kappa c} + e^{-6\kappa c}}{1 - e^{-6\kappa c}} =: G(c, \kappa).
\end{align*}

With $2\kappa c \approx 0.857615$ and $6\kappa c \approx 2.572845$,
\[
G(c, \kappa) \approx \frac{0.424180 + 0.076317}{1 - 0.076317} \approx 0.54185.
\]

We analyze $\text{RHS}(n) = (1 - 1/p_{n+1})^{-cp_n} - 1$.
\begin{align*}
\text{RHS}(n) &= \exp\left(cp_n \left(\frac{1}{p_{n+1}} + \frac{1}{2p_{n+1}^2} + O(p_{n+1}^{-3})\right)\right) - 1.
\end{align*}
We have $cp_n/p_{n+1} \geq a_0 c \approx 0.433096$. The leading term gives $\exp(0.433096) - 1 \approx 0.54204$.
The second-order term is small ($< 5.46 \times 10^{-7}$).

Thus, for $p_n \geq X_0$, $\text{RHS}(n) > 0.54203$.

The neglected contributions (tail terms and composite terms) are negligible compared to the gap of approximately $0.00018$ between the bounds $G(c,\kappa) \approx 0.54185$ and $\text{RHS}(n) > 0.54203$.

Comparing bounds, $T_n^*(cp_n) < 0.5419 < 0.54203 < \text{RHS}(n)$.
Since the inequality holds for all sufficiently large $n$, we conclude $C \leq 0.4332$.
\end{proof}

\begin{theorem}[Twin prime lower bound]
\label{thm:twin}
Under the twin prime conjecture, $C > \log \psi \approx 0.38225$, where $\psi$ denotes the supergolden ratio, the unique real root of $x^3 - x^2 - 1 = 0$.
\end{theorem}

\begin{proof}
By definition (Proposition \ref{prop:existence_sn}), $s_n^*$ satisfies $h(s_n^*) = p_{n+1} - 1$. Recalling $h(s) = p_{n+1}(1 + T_n^*(s))^{-1/s}$, we rearrange this equality at $s=s_n^*$ to obtain the exact equilibrium equation:
\[
1 + T_n^*(s_n^*) = \left(1 - \frac{1}{p_{n+1}}\right)^{-s_n^*}.
\]

Assuming the twin prime conjecture, consider the subsequence $\mathcal{T} = \{n : g_{n+1} = 2\}$. Define $C_{\text{twin}} := \limsup_{n \in \mathcal{T}} \sigma_n^*$. By PNT, $C_{\text{twin}} = \limsup_{n \in \mathcal{T}} s_n^*/p_{n+1}$.

We analyze the asymptotic behavior as $n \to \infty$ along a subsequence achieving the limsup.
\[
\left(1 - \frac{1}{p_{n+1}}\right)^{-s_n^*} = \exp\left(-s_n^* \log\left(1 - \frac{1}{p_{n+1}}\right)\right) = \exp\left(\frac{s_n^*}{p_{n+1}} + O\left(\frac{s_n^*}{p_{n+1}^2}\right)\right).
\]
Since $s_n^* = O(p_n)$ (Theorem \ref{thm:nagura}), the error term vanishes. Taking the limit, the equilibrium equation becomes
\[
\lim_{n\to\infty, n\in\mathcal{T}} (1 + T_n^*(s_n^*)) = \exp(C_{\text{twin}}).
\]

We analyze $T_n^*(s)$. In the limiting regime as $n \to \infty$ along the subsequence $\mathcal{T}$, composite contributions to $T_n^*(s_n^*)$ decay exponentially faster than prime contributions. The smallest composite coprime to $P_n$ exceeding $p_{n+1}$ is at least $p_{n+1}^2$. Since $s_n^* = O(p_n)$, the term $(p_{n+1}/k)^{s_n^*}$ for $k \geq p_{n+1}^2$ vanishes exponentially fast (e.g., $(1/p_{n+1})^{O(p_n)}$). Thus, only prime contributions remain in the limit.

$T_n^*(s) \geq \sum_{j=1}^{\infty} (1 + S_j/p_{n+1})^{-s}$. In the limit,
\[
T_n^*(s_n^*) \to H(C_{\text{twin}}) := \sum_{j=1}^{\infty} e^{-C_{\text{twin}} \cdot \overline{G}_j},
\]
where $\overline{G}_j$ are the limiting normalized cumulative gaps. The equilibrium is $\exp(C_{\text{twin}}) = 1 + H(C_{\text{twin}})$.

For twin primes ($g_{n+1} = 2$), admissibility modulo 3 forces $g_{n+2} \geq 4$. If $p_{n+1}, p_{n+1}+2$ are primes $>3$, then $p_{n+1} \equiv 2 \pmod 3$. Then $p_{n+1}+4 \equiv 0 \pmod 3$. Thus $S_1 = 2$ and $S_2 \geq 6$. The densest admissible pattern is $(2, 4, 2, 4, \ldots)$.

Therefore, $H(C) \geq e^{-2C} + e^{-6C} + e^{-8C} + \cdots > e^{-2C}$.

From $\exp(C) = 1 + H(C) > 1 + e^{-2C}$, we obtain $e^{3C} - e^{2C} - 1 > 0$.

Let $y = e^C > 1$. Define $f(y) = y^3 - y^2 - 1$. $f(y)$ is strictly increasing for $y > 2/3$. The unique real root is the supergolden ratio $\psi$. By Cardano's formula,
\[
\psi = \frac{1}{3}\left(1 + \sqrt[3]{\frac{29+3\sqrt{93}}{2}} + \sqrt[3]{\frac{29-3\sqrt{93}}{2}}\right) \approx 1.46557.
\]
Since $f(y) > 0$ requires $y > \psi$, we conclude $e^C > \psi$, giving $C > \log \psi \approx 0.38225$.
\end{proof}

\begin{remark}[The constant $C$ and prime constellations]
\label{rem:const_C}
The constant $C$ encodes which admissible prime constellations occur infinitely often. Heuristically, $C = \limsup \sigma_n$ is determined by the densest admissible $k$-tuple realized infinitely often. For twin primes (gap 2), the bound is $C \geq \log \psi \approx 0.3823$. Denser configurations (if they exist infinitely often) would yield larger values. Proving this rigorously requires resolving the Hardy--Littlewood $k$-tuple conjectures \cite{HardyLittlewood1923}, which predict the asymptotic density of each admissible pattern.
\end{remark}

\subsection{Empirical distribution of $\sigma_n$}

\begin{figure}[htbp]
\centering
\includegraphics[width=0.8\linewidth]{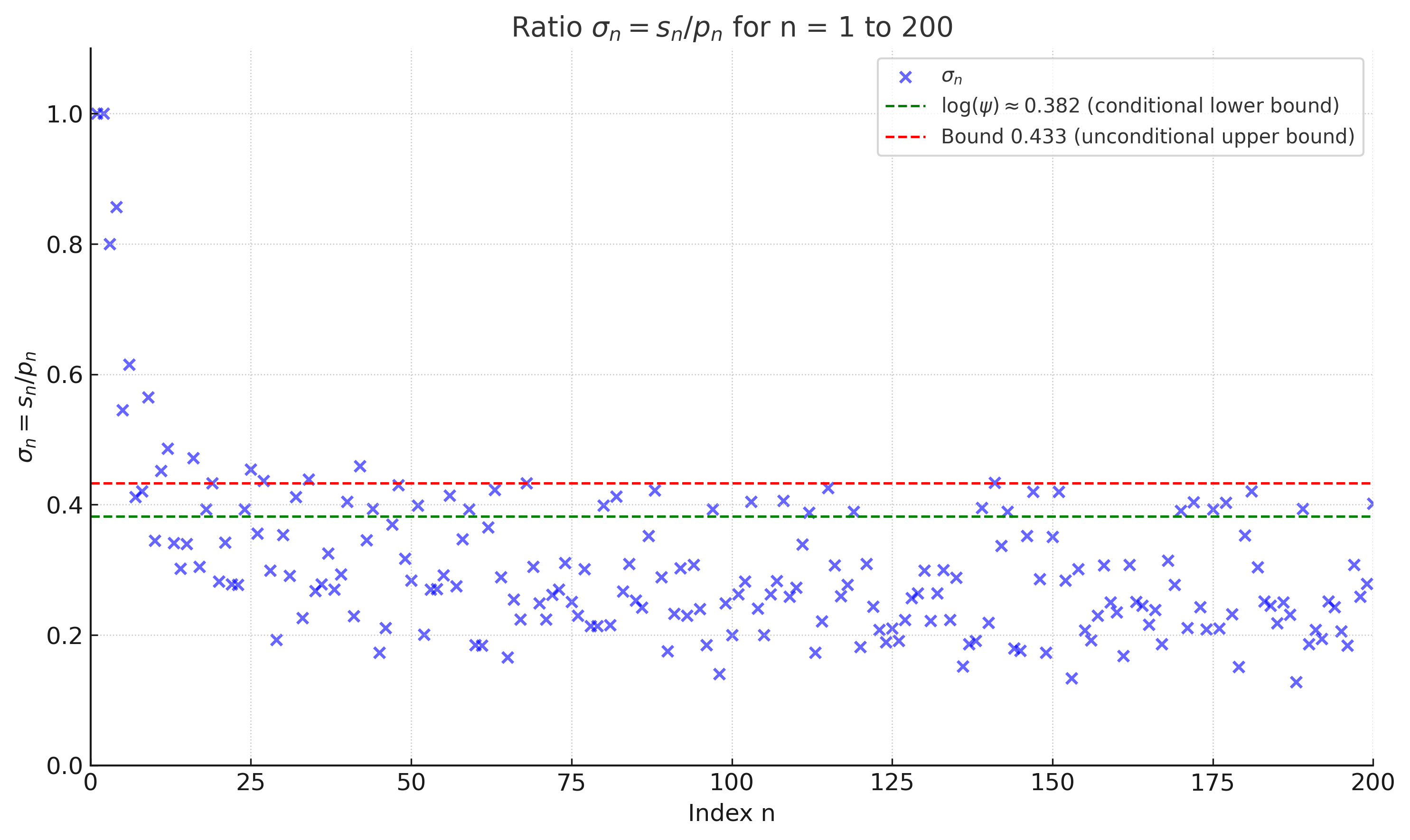}
\caption{The ratio $\sigma_n = s_n/p_n$ for $n=1$ to $200$. The horizontal lines indicate the theoretical bounds for $C = \limsup \sigma_n$: $\log \psi \approx 0.382$ (conditional lower bound) and $0.4332$ (unconditional upper bound).}
\label{fig:ratio}
\end{figure}

While the asymptotic analysis focuses on the extremal behavior ($\liminf$ and $\limsup$), we examine the empirical distribution of $\sigma_n$ for $n=1$ to 200 (Figure \ref{fig:ratio}) to understand its typical behavior.

\begin{figure}[htbp]
\centering
\includegraphics[width=0.8\linewidth]{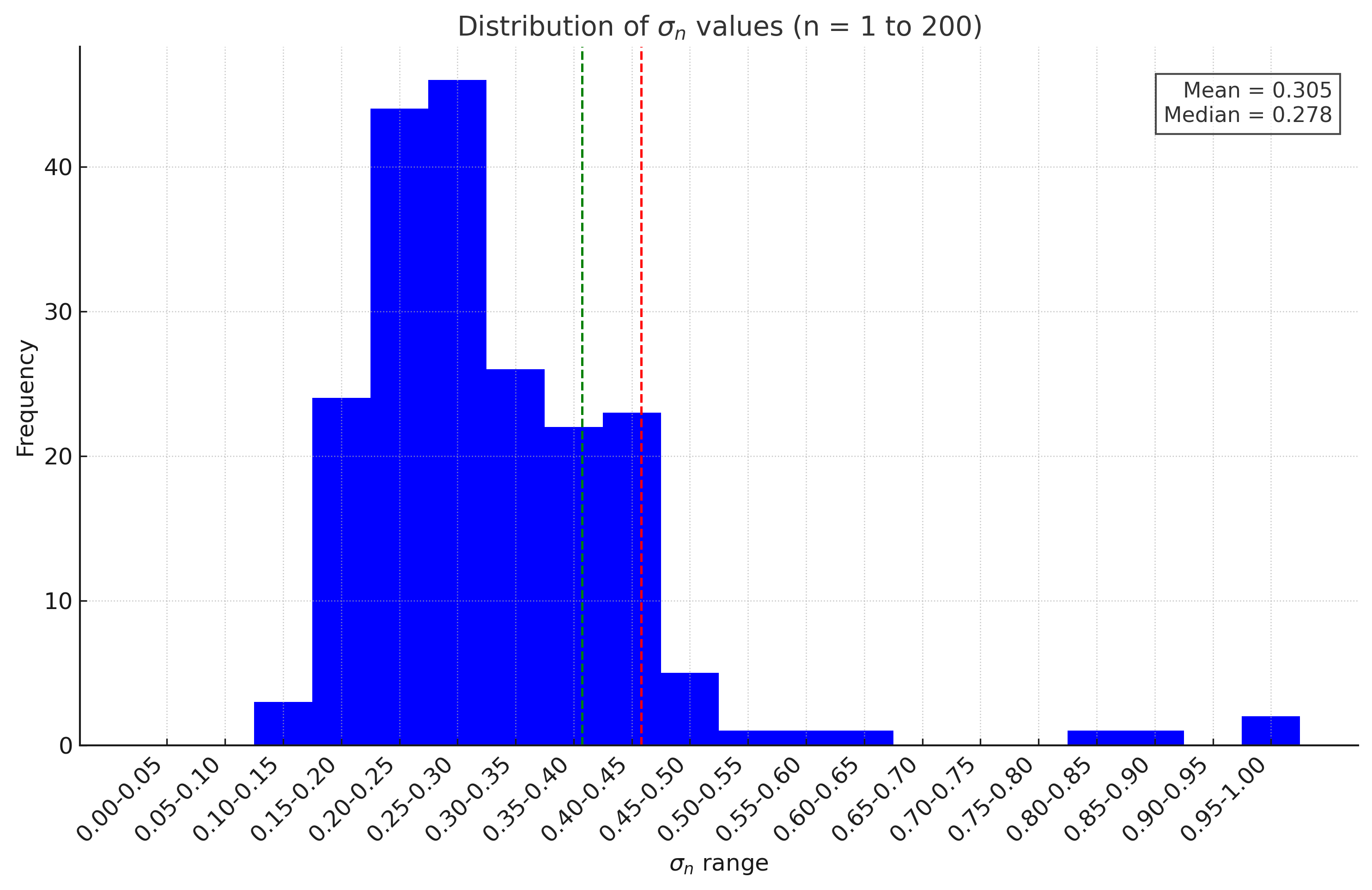}
\caption{Histogram of $\sigma_n$ for $n=1$ to 200. The distribution is heavily skewed right. Vertical dashed lines indicate the theoretical bounds for $C$: $\log \psi \approx 0.382$ (green) and $0.4332$ (red).}
\label{fig:histogram}
\end{figure}

The raw data histogram (Figure \ref{fig:histogram}) exhibits significant positive skewness (Mean $\approx 0.305$), largely due to the initial values.

To analyze the main body of the distribution, we remove outliers using the robust Tukey's fences method. The upper fence (Q3 + 1.5 $\times$ IQR) is calculated at $\approx 0.58$. This identifies 5 outliers corresponding to the initial values: $\sigma_n$ for $n=1, 2, 3, 4, 6$ (e.g., $\sigma_1=\sigma_2=1$). These early values occur because $p_n$ is very small ($p_n \leq 13$), and the asymptotic regime analyzed in Section 4.2 has not yet stabilized.

\begin{figure}[htbp]
\centering
\includegraphics[width=0.8\linewidth]{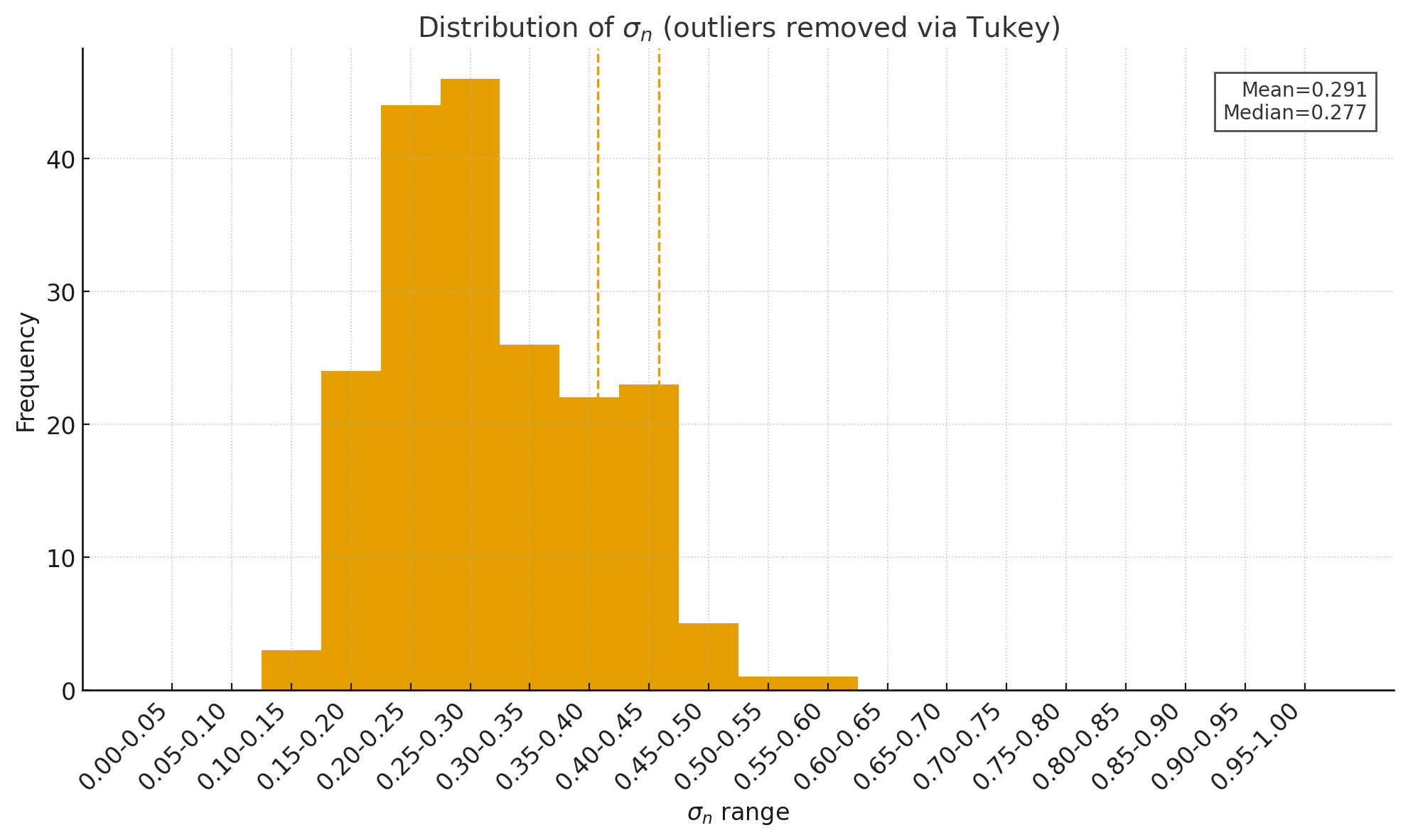}
\caption{Histogram of $\sigma_n$ ($n=1$ to 200) after removing 5 outliers identified by Tukey's fences. The core distribution is unimodal and moderately right-skewed, closely resembling a Beta(7.64, 18.62) distribution.}
\label{fig:histogram_trimmed}
\end{figure}

After removing these 5 points, the remaining 195 values (Figure \ref{fig:histogram_trimmed}) exhibit a clearer structure:
\begin{itemize}[leftmargin=*, nosep]
    \item Mean $\approx 0.291$, Median $\approx 0.277$, Standard Deviation $\approx 0.087$.
    \item Skewness $\approx +0.53$ (moderately right-skewed).
    \item Range [0.128, 0.565].
\end{itemize}

The sequence $\sigma_n$ is theoretically bounded within $[0, 1]$ (and asymptotically within $[0, C]$). The Beta distribution is a natural candidate for modeling this behavior due to its bounded support and flexible shape.

We fit a Beta distribution to the trimmed data using the method of moments. The estimated parameters are $\text{Beta}(\alpha \approx 7.64, \beta \approx 18.62)$. The theoretical mode of this distribution is $(\alpha-1)/(\alpha+\beta-2) \approx 0.274$, which aligns well with the empirical median and the modal bin [0.25, 0.30].

The "Beta-like" behavior suggests that while $\sigma_n$ can approach 0 (due to large gaps) and approach $C$ (due to dense constellations), the values predominantly cluster around a central tendency near 0.27--0.29. This centralization indicates a statistical regularity in the required precision $s_n$ relative to $p_n$, despite the underlying irregularity of prime gaps.

\section{Fixed coefficient formulas}

\subsection{The threshold constant $c_0$}

To determine when a fixed coefficient $c$ suffices, we analyze the worst-case scenario using the unconditional integral bound. For $s = cp_n$, the critical inequality asymptotically becomes
\[
e^{-cg_{n+1}} \left(1 + \frac{1}{c}\right) < e^c - 1.
\]
The worst case occurs for the smallest gap, $g_{n+1} = 2$. The threshold is defined by equality.

\begin{definition}[Critical threshold]
\label{def:c0}
Let $c_0$ be the unique positive solution to $e^{-2c}(1 + 1/c) = e^c - 1$.
Numerically, $c_0 \approx 0.5955790832$.
\end{definition}
\begin{proof}[Proof of Uniqueness]
Let $G(c) = e^c - 1 - e^{-2c}(1 + 1/c)$. We analyze the derivative for $c>0$:
\[
G'(c) = e^c + 2e^{-2c}\left(1 + \frac{1}{c}\right) + e^{-2c}\left(\frac{1}{c^2}\right).
\]
Since all terms are strictly positive for $c>0$, $G(c)$ is strictly increasing. As $\lim_{c\to 0^+} G(c) = -\infty$ and $\lim_{c\to \infty} G(c) = \infty$, there exists a unique root $c_0$.
\end{proof}

\begin{theorem}[Fixed coefficient]
\label{thm:fixed}
For any $c > c_0$, there exists $N_0(c)$ such that for all $n \geq N_0(c)$, the formula with $s=cp_n$ correctly computes $p_{n+1}$.
\end{theorem}

\begin{proof}
Let $s = cp_n$ with $c > c_0$. We established $\lim_{n\to\infty} \text{RHS}(n) = e^c - 1$ and $\limsup_{n\to\infty} T_n^*(cp_n) \leq e^{-2c}(1+1/c)$. By Definition \ref{def:c0}, the latter is strictly less than the former.

Let $\epsilon(c) = (e^c - 1) - e^{-2c}(1+1/c) > 0$. We control the error terms using effective bounds on prime gaps. By Baker--Harman--Pintz \cite{BakerHarmanPintz2001}, $g_n = O(p_n^{0.525})$, implying $p_n/p_{n+1} = 1 - O(p_n^{-0.475})$.

Using these effective bounds in the asymptotic expansions, the required inequality $T_n^*(s) < \text{RHS}(n)$ holds when the error terms (bounded by $O(p_n^{-0.475})$) are smaller than $\epsilon(c)$, which occurs for all sufficiently large $n$.
\end{proof}

\begin{conjecture}[General fixed coefficient]
\label{conj:general}
If $c > C = \limsup \sigma_n$, then the formula with $s = cp_n$ holds for all sufficiently large $n$.
\end{conjecture}
This conjecture is motivated by the definition of $C$ as the asymptotic maximum required normalized precision. Numerical verification for $c$ slightly above $C$ (e.g., $c=0.45$) supports this conjecture for the available data ($n \leq 200$).

\section{Extension to Dirichlet L-functions}

Haley \cite{Haley2013} extended Keller's formula to Dirichlet L-functions. We make this extension constructive and show that the sign of the L-series carries arithmetic information about $p_{n+1}$ accessible with reduced precision.

\subsection{Effective formula via L-functions}

\begin{definition}
The non-principal Dirichlet character modulo 4 is $\chi_4(n) = 1$ if $n \equiv 1 \pmod 4$, $-1$ if $n \equiv 3 \pmod 4$, and $0$ otherwise. Its L-function is $L(s, \chi_4) = 1 - 3^{-s} + 5^{-s} - \cdots$.
\end{definition}

\begin{theorem}[Effective L-function formula]
\label{thm:dirichlet_main}
Define
\[
V_n(s) = L(s, \chi_4) \prod_{k=1}^{n} (1 - \chi_4(p_k)p_k^{-s}).
\]
For all $n \geq 1$, taking $s = 2p_n$ yields
\[
p_{n+1} = \left\lceil \left|V_n(2p_n) - 1\right|^{-1/(2p_n)} \right\rceil.
\]
\end{theorem}

The proof parallels that of Theorem~\ref{thm:bertrand}, using the Euler product for $L(s,\chi_4)$ and integral estimates combined with Bertrand's postulate (as detailed in Theorem~\ref{thm:dirichlet_sign}).

\subsection{Arithmetic information from signs}

The Euler product expansion gives
\[
V_n(s) = 1 + \frac{\chi_4(p_{n+1})}{p_{n+1}^s} + \frac{\chi_4(p_{n+2})}{p_{n+2}^s} + \cdots
\]
The sign of $V_n(s) - 1$ reveals $p_{n+1} \pmod 4$ when the leading term dominates.

\begin{theorem}[Residue prediction from sign]
\label{thm:dirichlet_sign}
For $s = 2p_n$, the sign of $V_n(s) - 1$ determines $p_{n+1} \pmod 4$.
\end{theorem}

\begin{proof}
We need the leading term magnitude $1/p_{n+1}^s$ to dominate the tail. With $s = 2p_n$, we bound the tail using $|\chi_4(k)| \leq 1$ and Lemma \ref{lem:integral}(1):
\[
|\text{tail}| < \frac{(p_{n+2}-1)^{1-s}}{s-1}.
\]
We verify $1/p_{n+1}^s > (p_{n+2}-1)^{1-s}/(s-1)$, which rearranges to
\[
\frac{p_{n+2}-1}{s-1} \left(\frac{p_{n+1}}{p_{n+2}-1}\right)^s < 1.
\]
Bertrand's postulate gives $p_{n+2} < 4p_n$. The coefficient $(p_{n+2}-1)/(2p_n-1) < 2.2$ for $n \geq 2$. Since $p_{n+2} \geq p_{n+1} + 2$, the ratio term is bounded by $(1 - 1/(2p_n))^{2p_n} < e^{-1}$. The product is less than $2.2 \cdot e^{-1} \approx 0.81 < 1$.
The leading term dominates, so $\text{sign}(V_n(2p_n) - 1) = \text{sign}(\chi_4(p_{n+1}))$.
\end{proof}

\subsection{Minimal parameters for sign determination}

Define $s'_n$ as the minimal integer parameter for correct sign prediction, and $\sigma'_n = s'_n/p_n$.

\begin{theorem}[Lambert W lower bound]
\label{thm:lambert}
Let $C' := \limsup_{n \to \infty} \sigma'_n$. Under the twin prime conjecture,
\[
C' > W(1) \approx 0.567143,
\]
where $W$ is the Lambert W function. The first-order correction for finite $p_n$ along twin prime configurations is
\[
\sigma'_n \approx W(1) + \frac{K}{p_n}, \quad \text{where } K = \frac{1 + 1/W(1) + W(1)/2}{1 + 1/W(1)} \approx 1.1026235.
\]
\end{theorem}

\begin{proof}
The threshold $s^*$ occurs when the leading term equals the integral tail bound. Assuming the twin prime conjecture, consider the subsequence where $p_{n+2} = p_{n+1} + 2$. Let $p = p_{n+1}$ and $s^* = cp$. The equilibrium equation is $(s^*-1)(p+1)^{s^*-1} = p^{s^*}$.

Asymptotic expansion for large $p$ yields
\begin{equation}
\label{eq:logc_c_expansion}
\log c + c = \frac{1}{p}\left(1 + \frac{1}{c} + \frac{c}{2}\right) + O(p^{-2}).
\end{equation}
As $p \to \infty$, we get $\log c + c = 0$, or $c e^c = 1$. The solution is $c_0 = W(1)$, establishing $C' \geq W(1)$.

The first-order correction $K$ is derived by substituting $c = c_0 + K/p$ into the expansion (\ref{eq:logc_c_expansion}) and solving for $K$, utilizing the identity $1+c_0 = c_0(1+1/c_0)$. This yields the explicit expression for $K$. Since $K>0$, we have $C' > W(1)$.
\end{proof}

\section{Discussion and Computational Aspects}

\subsection{Prime constellations and jumping champions}

Odlyzko, Rubinstein, and Wolf \cite{OdlyzkoRubinstein1999} studied jumping champions (most frequent prime gaps up to $x$). Our constant $C$ is influenced by gap distribution differently: jumping champions describe local frequency, while $C$ is determined by configurations producing the largest $\sigma_n$ infinitely often. Understanding how $C$ relates quantitatively to the transition points of jumping champions remains an intriguing open problem, requiring deep analysis of how local gap frequencies (governed by the Hardy--Littlewood singular series) connect to the global extremal behavior captured by $\limsup \sigma_n$.

\subsection{Computational complexity}

We analyze the complexity of computing $p_{n+1}$ using the recurrence, which requires finding the minimal $s_n$ by iterating $s$ from 1 up to $p_n$ (worst case). The cost is dominated by the evaluation of the Euler product. With multi-precision arithmetic at $O(s \log s)$ bits (since $s = O(p_n)$ and $p_n \sim n \log n$), each of the $n$ primorial factor computations (exponentiation and multiplication) requires $O(s^2 \log s)$ bit operations (assuming standard multiplication complexity). The cost per iteration $s$ is thus $O(n s^2 \log s)$. Summing over the iterations, the total complexity is $\sum_{s=1}^{p_n} O(n s^2 \log s) = O(n p_n^3 \log p_n)$ bit operations. While polynomial in the input size, it remains impractical compared to sieving methods which require $O(p_n \log \log p_n)$ operations. The significance of the formula is theoretical, introducing the sequence $(s_n)$ which possesses intrinsic number-theoretic interest.

\subsection{Conclusion}

We have transformed the asymptotic Golomb--Keller formula into a constructive recurrence by establishing the explicit bound $s_n \leq p_n$. The analysis of $\sigma_n$ uncovers a deep connection to prime constellations, with its limit superior $C$ bounded by $0.3823 \lesssim C \leq 0.4332$. The appearance of the supergolden ratio in the conditional lower bound is novel in this context. The empirical distribution of $\sigma_n$ suggests a Beta-like behavior, indicating centralization despite the irregularity of prime gaps. The method extends to Dirichlet L-functions, yielding a practical tool for predicting prime residues.

\section*{Acknowledgments}

The author thanks Philippe Caldero for initial inspiration.

\appendix
\section{Numerical Data}

The following table presents the values of $s_n$ and $\sigma_n$ for $n=1$ to 200, computed using PARI/GP \cite{PARI2023} with 100 digits precision. This data verifies Theorem~\ref{thm:nagura} ($s_n \leq p_n$) and forms the basis for the empirical analysis in Section 4.3.

\begin{center}
\begin{small}
\begin{longtable}{@{}crrrccc@{}}
\caption{Values of $s_n$ and $\sigma_n$ for $n=1$ to 200.} \label{tab:data} \\
\toprule
$n$ & $p_n$ & $p_{n+1}$ & $g_n$ & $g_{n+1}$ & $s_n$ & $\sigma_n = s_n/p_n$ \\
\midrule
\endfirsthead
\toprule
$n$ & $p_n$ & $p_{n+1}$ & $g_n$ & $g_{n+1}$ & $s_n$ & $\sigma_n = s_n/p_n$ \\
\midrule
\endhead
\bottomrule
\endfoot
1 & 2 & 3 & 1 & 2 & 2 & 1.000 \\
2 & 3 & 5 & 2 & 2 & 3 & 1.000 \\
3 & 5 & 7 & 2 & 4 & 4 & 0.800 \\
4 & 7 & 11 & 4 & 2 & 6 & 0.857 \\
5 & 11 & 13 & 2 & 4 & 6 & 0.545 \\
6 & 13 & 17 & 4 & 2 & 8 & 0.615 \\
7 & 17 & 19 & 2 & 4 & 7 & 0.412 \\
8 & 19 & 23 & 4 & 6 & 8 & 0.421 \\
9 & 23 & 29 & 6 & 2 & 13 & 0.565 \\
10 & 29 & 31 & 2 & 6 & 10 & 0.345 \\
11 & 31 & 37 & 6 & 4 & 14 & 0.452 \\
12 & 37 & 41 & 4 & 2 & 18 & 0.486 \\
13 & 41 & 43 & 2 & 4 & 14 & 0.341 \\
14 & 43 & 47 & 4 & 6 & 13 & 0.302 \\
15 & 47 & 53 & 6 & 6 & 16 & 0.340 \\
16 & 53 & 59 & 6 & 2 & 25 & 0.472 \\
17 & 59 & 61 & 2 & 6 & 18 & 0.305 \\
18 & 61 & 67 & 6 & 4 & 24 & 0.393 \\
19 & 67 & 71 & 4 & 2 & 29 & 0.433 \\
20 & 71 & 73 & 2 & 6 & 20 & 0.282 \\
21 & 73 & 79 & 6 & 4 & 25 & 0.342 \\
22 & 79 & 83 & 4 & 6 & 22 & 0.278 \\
23 & 83 & 89 & 6 & 8 & 23 & 0.277 \\
24 & 89 & 97 & 8 & 4 & 35 & 0.393 \\
25 & 97 & 101 & 4 & 2 & 44 & 0.454 \\
26 & 101 & 103 & 2 & 4 & 36 & 0.356 \\
27 & 103 & 107 & 4 & 2 & 45 & 0.437 \\
28 & 107 & 109 & 2 & 4 & 32 & 0.299 \\
29 & 109 & 113 & 4 & 14 & 21 & 0.193 \\
30 & 113 & 127 & 14 & 4 & 40 & 0.354 \\
31 & 127 & 131 & 4 & 6 & 37 & 0.291 \\
32 & 131 & 137 & 6 & 2 & 54 & 0.412 \\
33 & 137 & 139 & 2 & 10 & 31 & 0.226 \\
34 & 139 & 149 & 10 & 2 & 61 & 0.439 \\
35 & 149 & 151 & 2 & 6 & 40 & 0.268 \\
36 & 151 & 157 & 6 & 6 & 42 & 0.278 \\
37 & 157 & 163 & 6 & 4 & 51 & 0.325 \\
38 & 163 & 167 & 4 & 6 & 44 & 0.270 \\
39 & 167 & 173 & 6 & 6 & 49 & 0.293 \\
40 & 173 & 179 & 6 & 2 & 70 & 0.405 \\
41 & 179 & 181 & 2 & 10 & 41 & 0.229 \\
42 & 181 & 191 & 10 & 2 & 83 & 0.459 \\
43 & 191 & 193 & 2 & 4 & 66 & 0.346 \\
44 & 193 & 197 & 4 & 2 & 76 & 0.394 \\
45 & 197 & 199 & 2 & 12 & 34 & 0.173 \\
46 & 199 & 211 & 12 & 12 & 42 & 0.211 \\
47 & 211 & 223 & 12 & 4 & 78 & 0.370 \\
48 & 223 & 227 & 4 & 2 & 96 & 0.430 \\
49 & 227 & 229 & 2 & 4 & 72 & 0.317 \\
50 & 229 & 233 & 4 & 6 & 65 & 0.284 \\
51 & 233 & 239 & 6 & 2 & 93 & 0.399 \\
52 & 239 & 241 & 2 & 10 & 48 & 0.201 \\
53 & 241 & 251 & 10 & 6 & 65 & 0.270 \\
54 & 251 & 257 & 6 & 6 & 68 & 0.271 \\
55 & 257 & 263 & 6 & 6 & 75 & 0.292 \\
56 & 263 & 269 & 6 & 2 & 109 & 0.414 \\
57 & 269 & 271 & 2 & 6 & 74 & 0.275 \\
58 & 271 & 277 & 6 & 4 & 94 & 0.347 \\
59 & 277 & 281 & 4 & 2 & 109 & 0.393 \\
60 & 281 & 283 & 2 & 10 & 52 & 0.185 \\
61 & 283 & 293 & 10 & 14 & 52 & 0.184 \\
62 & 293 & 307 & 14 & 4 & 107 & 0.365 \\
63 & 307 & 311 & 4 & 2 & 130 & 0.423 \\
64 & 311 & 313 & 2 & 4 & 90 & 0.289 \\
65 & 313 & 317 & 4 & 14 & 52 & 0.166 \\
66 & 317 & 331 & 14 & 6 & 81 & 0.255 \\
67 & 331 & 337 & 6 & 10 & 74 & 0.224 \\
68 & 337 & 347 & 10 & 2 & 146 & 0.433 \\
69 & 347 & 349 & 2 & 4 & 106 & 0.305 \\
70 & 349 & 353 & 4 & 6 & 87 & 0.249 \\
71 & 353 & 359 & 6 & 8 & 79 & 0.224 \\
72 & 359 & 367 & 8 & 6 & 94 & 0.262 \\
73 & 367 & 373 & 6 & 6 & 99 & 0.270 \\
74 & 373 & 379 & 6 & 4 & 116 & 0.311 \\
75 & 379 & 383 & 4 & 6 & 95 & 0.251 \\
76 & 383 & 389 & 6 & 8 & 88 & 0.230 \\
77 & 389 & 397 & 8 & 4 & 117 & 0.301 \\
78 & 397 & 401 & 4 & 8 & 85 & 0.214 \\
79 & 401 & 409 & 8 & 10 & 86 & 0.214 \\
80 & 409 & 419 & 10 & 2 & 163 & 0.399 \\
81 & 419 & 421 & 2 & 10 & 90 & 0.215 \\
82 & 421 & 431 & 10 & 2 & 174 & 0.413 \\
83 & 431 & 433 & 2 & 6 & 115 & 0.267 \\
84 & 433 & 439 & 6 & 4 & 134 & 0.309 \\
85 & 439 & 443 & 4 & 6 & 111 & 0.253 \\
86 & 443 & 449 & 6 & 8 & 107 & 0.242 \\
87 & 449 & 457 & 8 & 4 & 158 & 0.352 \\
88 & 457 & 461 & 4 & 2 & 193 & 0.422 \\
89 & 461 & 463 & 2 & 4 & 133 & 0.289 \\
90 & 463 & 467 & 4 & 12 & 81 & 0.175 \\
91 & 467 & 479 & 12 & 8 & 109 & 0.233 \\
92 & 479 & 487 & 8 & 4 & 145 & 0.303 \\
93 & 487 & 491 & 4 & 8 & 112 & 0.230 \\
94 & 491 & 499 & 8 & 4 & 151 & 0.308 \\
95 & 499 & 503 & 4 & 6 & 120 & 0.240 \\
96 & 503 & 509 & 6 & 12 & 93 & 0.185 \\
97 & 509 & 521 & 12 & 2 & 200 & 0.393 \\
98 & 521 & 523 & 2 & 18 & 73 & 0.140 \\
99 & 523 & 541 & 18 & 6 & 130 & 0.249 \\
100 & 541 & 547 & 6 & 10 & 108 & 0.200 \\
101 & 547 & 557 & 10 & 6 & 144 & 0.263 \\
102 & 557 & 563 & 6 & 6 & 157 & 0.282 \\
103 & 563 & 569 & 6 & 2 & 228 & 0.405 \\
104 & 569 & 571 & 2 & 6 & 137 & 0.241 \\
105 & 571 & 577 & 6 & 10 & 114 & 0.200 \\
106 & 577 & 587 & 10 & 6 & 152 & 0.263 \\
107 & 587 & 593 & 6 & 6 & 166 & 0.283 \\
108 & 593 & 599 & 6 & 2 & 241 & 0.406 \\
109 & 599 & 601 & 2 & 6 & 155 & 0.259 \\
110 & 601 & 607 & 6 & 6 & 164 & 0.273 \\
111 & 607 & 613 & 6 & 4 & 206 & 0.339 \\
112 & 613 & 617 & 4 & 2 & 238 & 0.388 \\
113 & 617 & 619 & 2 & 12 & 107 & 0.173 \\
114 & 619 & 631 & 12 & 10 & 137 & 0.221 \\
115 & 631 & 641 & 10 & 2 & 269 & 0.426 \\
116 & 641 & 643 & 2 & 4 & 197 & 0.307 \\
117 & 643 & 647 & 4 & 6 & 167 & 0.260 \\
118 & 647 & 653 & 6 & 6 & 179 & 0.277 \\
119 & 653 & 659 & 6 & 2 & 254 & 0.389 \\
120 & 659 & 661 & 2 & 12 & 120 & 0.182 \\
121 & 661 & 673 & 12 & 4 & 204 & 0.309 \\
122 & 673 & 677 & 4 & 6 & 164 & 0.244 \\
123 & 677 & 683 & 6 & 8 & 141 & 0.208 \\
124 & 683 & 691 & 8 & 10 & 129 & 0.189 \\
125 & 691 & 701 & 10 & 8 & 145 & 0.210 \\
126 & 701 & 709 & 8 & 10 & 134 & 0.191 \\
127 & 709 & 719 & 10 & 8 & 158 & 0.223 \\
128 & 719 & 727 & 8 & 6 & 185 & 0.257 \\
129 & 727 & 733 & 6 & 6 & 192 & 0.264 \\
130 & 733 & 739 & 6 & 4 & 219 & 0.299 \\
131 & 739 & 743 & 4 & 8 & 164 & 0.222 \\
132 & 743 & 751 & 8 & 6 & 196 & 0.264 \\
133 & 751 & 757 & 6 & 4 & 225 & 0.300 \\
134 & 757 & 761 & 4 & 8 & 169 & 0.223 \\
135 & 761 & 769 & 8 & 4 & 219 & 0.288 \\
136 & 769 & 773 & 4 & 14 & 117 & 0.152 \\
137 & 773 & 787 & 14 & 10 & 144 & 0.186 \\
138 & 787 & 797 & 10 & 12 & 150 & 0.191 \\
139 & 797 & 809 & 12 & 2 & 315 & 0.395 \\
140 & 809 & 811 & 2 & 10 & 177 & 0.219 \\
141 & 811 & 821 & 10 & 2 & 352 & 0.434 \\
142 & 821 & 823 & 2 & 4 & 277 & 0.337 \\
143 & 823 & 827 & 4 & 2 & 320 & 0.389 \\
144 & 827 & 829 & 2 & 10 & 149 & 0.180 \\
145 & 829 & 839 & 10 & 14 & 146 & 0.176 \\
146 & 839 & 853 & 14 & 4 & 295 & 0.352 \\
147 & 853 & 857 & 4 & 2 & 358 & 0.420 \\
148 & 857 & 859 & 2 & 4 & 245 & 0.286 \\
149 & 859 & 863 & 4 & 14 & 149 & 0.173 \\
150 & 863 & 877 & 14 & 4 & 303 & 0.351 \\
151 & 877 & 881 & 4 & 2 & 368 & 0.420 \\
152 & 881 & 883 & 2 & 4 & 250 & 0.284 \\
153 & 883 & 887 & 4 & 20 & 118 & 0.134 \\
154 & 887 & 907 & 20 & 4 & 267 & 0.301 \\
155 & 907 & 911 & 4 & 8 & 188 & 0.207 \\
156 & 911 & 919 & 8 & 10 & 175 & 0.192 \\
157 & 919 & 929 & 10 & 8 & 211 & 0.230 \\
158 & 929 & 937 & 8 & 4 & 285 & 0.307 \\
159 & 937 & 941 & 4 & 6 & 234 & 0.250 \\
160 & 941 & 947 & 6 & 6 & 221 & 0.235 \\
161 & 947 & 953 & 6 & 14 & 159 & 0.168 \\
162 & 953 & 967 & 14 & 4 & 294 & 0.308 \\
163 & 967 & 971 & 4 & 6 & 243 & 0.251 \\
164 & 971 & 977 & 6 & 6 & 238 & 0.245 \\
165 & 977 & 983 & 6 & 8 & 211 & 0.216 \\
166 & 983 & 991 & 8 & 6 & 235 & 0.239 \\
167 & 991 & 997 & 6 & 12 & 184 & 0.186 \\
168 & 997 & 1009 & 12 & 4 & 313 & 0.314 \\
169 & 1009 & 1013 & 4 & 6 & 279 & 0.277 \\
170 & 1013 & 1019 & 6 & 2 & 396 & 0.391 \\
171 & 1019 & 1021 & 2 & 10 & 215 & 0.211 \\
172 & 1021 & 1031 & 10 & 2 & 413 & 0.404 \\
173 & 1031 & 1033 & 2 & 6 & 250 & 0.243 \\
174 & 1033 & 1039 & 6 & 10 & 216 & 0.209 \\
175 & 1039 & 1049 & 10 & 2 & 408 & 0.393 \\
176 & 1049 & 1051 & 2 & 10 & 220 & 0.210 \\
177 & 1051 & 1061 & 10 & 2 & 424 & 0.403 \\
178 & 1061 & 1063 & 2 & 6 & 246 & 0.232 \\
179 & 1063 & 1069 & 6 & 18 & 161 & 0.151 \\
180 & 1069 & 1087 & 18 & 4 & 377 & 0.353 \\
181 & 1087 & 1091 & 4 & 2 & 458 & 0.421 \\
182 & 1091 & 1093 & 2 & 4 & 332 & 0.304 \\
183 & 1093 & 1097 & 4 & 6 & 275 & 0.252 \\
184 & 1097 & 1103 & 6 & 6 & 269 & 0.245 \\
185 & 1103 & 1109 & 6 & 8 & 240 & 0.218 \\
186 & 1109 & 1117 & 8 & 6 & 277 & 0.250 \\
187 & 1117 & 1123 & 6 & 6 & 258 & 0.231 \\
188 & 1123 & 1129 & 6 & 22 & 144 & 0.128 \\
189 & 1129 & 1151 & 22 & 2 & 445 & 0.394 \\
190 & 1151 & 1153 & 2 & 10 & 214 & 0.186 \\
191 & 1153 & 1163 & 10 & 8 & 240 & 0.208 \\
192 & 1163 & 1171 & 8 & 10 & 226 & 0.194 \\
193 & 1171 & 1181 & 10 & 6 & 295 & 0.252 \\
194 & 1181 & 1187 & 6 & 6 & 287 & 0.243 \\
195 & 1187 & 1193 & 6 & 8 & 244 & 0.206 \\
196 & 1193 & 1201 & 8 & 12 & 219 & 0.184 \\
197 & 1201 & 1213 & 12 & 4 & 370 & 0.308 \\
198 & 1213 & 1217 & 4 & 6 & 314 & 0.259 \\
199 & 1217 & 1223 & 6 & 6 & 340 & 0.279 \\
200 & 1223 & 1229 & 6 & 2 & 492 & 0.402 \\
\end{longtable}
\end{small}
\end{center}

\end{document}